\documentclass[11pt,reqno,twoside]{amsart}

\usepackage[english]{babel} 
\usepackage{esint}
\usepackage[T1]{fontenc} 
\usepackage{amsmath}
\usepackage{amssymb} 
\usepackage{float}
\usepackage{amsfonts}
\usepackage{mathtools}
\usepackage{xcolor}
\usepackage[utf8]{inputenc}
\usepackage[mathcal]{eucal} 
\usepackage{enumerate}
\usepackage{amsthm} 
\usepackage{hyperref}
\usepackage[totalwidth=14cm,totalheight=20cm, hmarginratio=1:1]{geometry}

\newcommand{\R}{\mathbb{R}}

\newcommand{\SO}{\mathrm{SO}}
\newcommand{\h}{\mathbb{H}}
\newcommand{\PSL}{\mathbb{P}\mathrm{SL}}
\newcommand{\Sp}{\mathrm{Sp}}

\newcommand{\Pp}{\mathbb{P}}

\newcommand{\Mix}{\mathrm{Mix}}
\newcommand{\Flat}{\mathrm{Flat}}
\newcommand{\Area}{\mathrm{Area}}

\newtheorem{defi}{Definition}[section] 
\newtheorem{thm}[defi]{Theorem}
\newtheorem{cor}[defi]{Corollary} 
\newtheorem{lemma}[defi]{Lemma}

\newtheorem{prop}[defi]{Proposition}

\newtheorem{bigthm}{Theorem}

\theoremstyle{definition}
\newtheorem{rmk}{Remark}[section]

\DeclareMathAlphabet{\mathpzc}{OT1}{pzc}{m}{it}

\begin{document}

\title[Boundary of the Gothen components]{Boundary of the Gothen components}
\author{Charles Ouyang and Andrea Tamburelli}

\begin{abstract} In this short note we describe an interesting new phenomenon about the $\Sp(4,\R)$-character variety. Precisely, we show that the Hitchin component and all Gothen components share the same boundary in our length spectrum compactification.
\end{abstract}

\maketitle
\setcounter{tocdepth}{1}
\tableofcontents

\section*{Introduction}
This short note describes an interesting new phenomenon about the length spectrum compactification of $\Sp(4,\R)$-maximal representations of the fundamental group of a closed surface $S$. These were first introduced by Burger-Iozzi-Wienhard (\cite{BIW_maximal}) as generalizations of discrete and faithful representations into $\PSL(2,\R)$ and have been studied since then by many authors for their geometric and dynamical properties (\cite{BGPG_Hermitian}, \cite{BP_maximal}, \cite{maximal_Fenchel_Nielsen}). 
In particular, Gothen (\cite{gothen2001components}) showed that there are precisely $3\cdot 2^{2g}+2g-4$ connected components of maximal representations, of which only $2^{2g}+2g-3$ are smooth. Of these, $2^{2g}$ are isomorphic copies of the $\Sp(4,\R)$-Hitchin component, which were analyzed in a previous work (\cite{OT_Sp4}). The remaining $2g-3$ are referred to as Gothen components in the literature and are the main object of study of this paper. \\

Collier (\cite{collier2016maximal}), extending previous work of Labourie (\cite{Labourie_cyclic}), showed that for every maximal representation in $\Sp(4,\R)$ there is a unique equivariant conformal harmonic map from the universal cover of $S$ into the symmetric space $\Sp(4,\R)/\mathrm{U}(2)$. This allows for a parameterization of the connected components of maximal representations as bundles over the Teichm\"uller space of $S$ (\cite{AC_Sp4}). Exploiting the exceptional isomorphism $\Pp\Sp(4,\R)\cong \SO_{0}(2,3)$,  Collier-Tholozan-Toulisse (\cite{BTT}) interpreted the unique equivariant minimal immersion into the symmetric space as the Gauss map of the unique equivariant maximal surface (i.e with vanishing mean curvature) in the pseudo-hyperbolic space $\h^{2,2}$. Moreover, the different connected components of maximal representations correspond in this language to the different isomorphism classes that the orthogonal bundle to the maximal surface can have. \\

In a previous paper (\cite{OT_Sp4}), we used this pseudo-Riemannian point of view to define the length spectrum of an $\Sp(4,\R)$-Hitchin representation and describe the boundary of this component as projectivized mixed structures, that is hybrid geometric objects on the surface that are measured laminations on some collection of incompressible subsurfaces and flat metrics of finite area induced by meromorphic quartic differentials on the complement. In this note, we show that all Gothen components share the same boundary, precisely

\vspace{0.1cm}
\begin{bigthm}\label{thmA} For all $1\leq k \leq 2g-3$, the boundary of the Gothen component $\mathpzc{G}(k)$ can be identified with the space of projectivized mixed structures $\Pp\Mix_{4}(S)$.
\end{bigthm}
\vspace{0.1cm}

This behavior is in stark contrast to a result of Wolff (\cite{wolff2011connected}), where the projective measured laminations intersect as a closed nowhere dense set in the boundary of the non-Teichm{\"u}ller components of the $\PSL(2, \R)$-character variety. Moreover, Theorem \ref{thmA} is germane only to $\Sp(4, \R)$, as these exceptional Gothen components no longer appear in the $\Sp(2n, \R)$-character variety once $n$ is at least 3 (\cite{garcia2013higgs}). Garc{\'i}a-Prada-Gothen-Mundet showed there are only $3\cdot 2^{2g}$ maximal components in the $\Sp(2n, \R)$-character variety once $n$ exceeds 2.\\

The main idea behind the proof of Theorem \ref{thmA} lies in a comparison between the induced metric on the maximal surface associated to a maximal representation and the flat metric $|q|^{\frac{1}{2}}$ with cone singularities, where $q$ is the holomorphic quartic differential in the Hitchin base. We show that for any diverging sequence in $\mathpzc{G}(k)$, the length spectra of these two metrics (after a suitable renormalization) converge to that of mixed structures which enjoy the same decomposition into subsurfaces and coincide on their non-laminar part.  \\

As an application of our estimates, we also describe the length spectrum for the induced metrics on the unique minimal surfaces in $\Sp(4,\R)/\mathrm{U}(2)$.

\begin{bigthm}\label{thmB} Let $\rho_{n}$ be a sequence of representations in a Gothen component that leaves every compact set. Let $g_{n}$ be the corresponding sequence of induced metrics on the unique equivariant minimal surfaces in $\Sp(4,\R)/\mathrm{U}(2)$. Then the marked length spectrum of $g_{n}$ converges projectively to that of a mixed structure $\eta \in\mathbb{P}\mathrm{Mix}_{4}(S)$. Moreover, all such mixed structures can be attained in the limit.
\end{bigthm}

Similar results have been proven for the Blaschke metrics on affine spheres (\cite{OT}) and equivariant minimal Lagrangian surfaces in $\mathbb{H}^{2}\times \mathbb{H}^2$ (\cite{Charles_dPSL}), thus this work concludes the study of the length spectrum compactification of higher Teichm\"uller components for the classical Lie groups of rank $2$.

\subsection*{Acknowledgement} This paper was written while the first author was visiting Rice University during the academic year 2020-21. He thanks the Department of Mathematics for their kindness and hospitality. The second author acknowledges support from the National Science Foundation through grant DMS-2005501.

\section{Background material}

\subsection{Maximal representations}
Let $G$ be a Lie group, whose associated symmetric space is Hermitian. Then to any representation $\rho: \pi_{1}(S) \to G$, there is a smooth $\pi_{1}(S)$-equivariant map $\widetilde{f}_{\rho}: \widetilde{S} \to G/K$ given by taking any section of the flat bundle $E_{\rho}=\widetilde{S} \times_{\rho} G/K \rightarrow S$. The $G$-invariant two-form $\omega$ on $G/K$ pulls back via $\widetilde{f}_{\rho}$ to a well-defined two form on $S$. The \textit{Toledo invariant} of $\rho$ is defined as $T_{\rho} := \int_{S} \widetilde{f}_{\rho}^{*} \omega$. As the fiber of the bundle $E_{\rho}$ is contractible, a choice of any other section would yield another map differing by a $\rho$-equivariant homotopy, so that the integral is well-defined. In particular, one has the following Milnor-Wood type inequality: $|T_{\rho}| \leq (2g-2) \, \text{rank}(G)$ (see \cite{burger2005maximal}). As the Toledo invariant is constant on connected components of the representation variety, those attaining the maximal Toledo number are called \textit{maximal representations}. These have been extensively studied: see for instance \cite{BIW_maximal}, \cite{burger2005maximal}. For example, in the setting of $G=\text{PSL}(2,\mathbb{R})$, Goldman (\cite{goldmanthesis}) shows that each component is characterized by its Toledo number, with the two maximal components each being a copy of Teichm{\"u}ller space. Hence maximal components can be viewed as higher rank analogues of the classical Teichm{\"u}ller space.

When $G= \Sp(4,\R)$, Gothen (\cite{gothen2001components}) showed there are precisely 3 $\cdot 2^{2g} + 2g-4$ maximal components, of which only $2^{2g} + 2g-3$ are smooth (\cite{bradlow2012deformations}). Of these, $2^{2g}$ are isomorphic copies of the $\Sp(4, \R)$-Hitchin component, so the remaining $2g-3$ are often called Gothen components in the literature. Like their Hitchin counterparts, each representation $\rho$ in the Gothen components has a unique Riemann surface structure associated to them (see \cite{collier2016maximal}) so that the $\rho$-equivariant harmonic map from the universal cover of that Riemann surface into the symmetric space is conformal. Using the preferred Riemann surface, under the non-abelian Hodge correspondence, Gothen representations are associated to Higgs bundles of the form
$$ \mathcal{E} = N \oplus N K^{-1} \oplus N^{-1}K \oplus N^{-1}, \qquad \phi= \begin{pmatrix}
0 & 0 & 0 & \nu\\
1 & 0 & 0 & 0 \\
0 & \mu & 0 & 0\\
0 & 0 & 1 & 0 \\
\end{pmatrix}, $$
where $N$ is a line bundle for which $g-1< \deg N \leq 3g-3$ and $0 \neq \mu \in H^{0}(X, N^{-2}K^{3})$ and $\nu \in H^{0}(X, N^{2}K)$.
In particular, $\mu \nu \in H^{0}(X, K^{4})$. When $N = K^{3/2}$ and $\mu =1$, one recovers cyclic Higgs bundles corresponding to Hitchin $\Sp(4, \R)$-representations (see \cite{Labourie_cyclic}). Using the preferred Riemann surface, the Hitchin equations for the Higgs bundles corresponding to Gothen representations are given by

\begin{equation}\label{eq:Hitchin}
\begin{cases}
\Delta \psi_{1} =e^{\psi_{1} - \psi_{2}}-|\nu|^{2}e^{-2\psi_{1}}   \\ 
\Delta \psi_{2} =|\mu|^{2}e^{2\psi_{2}}- e^{\psi_{1} - \psi_{2}} 
\end{cases} \ ,
\end{equation}
where $e^{\psi_{1}}$ and $e^{\psi_{2}}$ are the local expressions of the hermitian harmonic metrics on the line bundles $N^{-1}$ and $N^{-1}K$ respectively. 

We remark that the family of Higgs bundles shown above do not always give rise to different representations. In fact, when $\deg N < 3g-3$ there is a natural $\mathbb{C}^{*}$ action that sends the pair $(\mu, \nu)$ to $(\lambda \mu, \lambda^{-1}\nu)$. Then, by a recent result of Alessandrini and Collier (\cite{AC_Sp4}), the Gothen component $\mathpzc{G}(k)$ is homeomorphic to the bundle over the  Teichm\"uller space of $S$ whose fiber is
\[
    \{ (\mu, \nu) \in H^{0}(X,N^{-2}K^{2}) \times H^{0}(X,N^{2}K)\} / \mathbb{C}^{*}
\]
where $N$ is a line bundle of degree $k$. 

Recall that the Hitchin fibration associates to any Higgs bundle the coefficients of the characteristic polynomial of the Higgs field, which are holomorphic differentials on $X$. In this case, the only non-zero differential is the holomorphic quartic differential $q=\mu \nu \in H^{0}(X, K^{4})$. Higgs bundles with the same quartic differentials are said to be in the same Hitchin fiber. Among these, the Higgs bundle
\[
    \mathcal{E} = K^{\frac{3}{2}} \oplus K^{\frac{1}{2}} \oplus K^{-\frac{1}{2}}K \oplus K^{-\frac{3}{2}}, \qquad \phi= \begin{pmatrix}
0 & 0 & 0 & q\\
1 & 0 & 0 & 0 \\
0 & 1 & 0 & 0\\
0 & 0 & 1 & 0 \\
\end{pmatrix}
\]
is the unique one whose associated representation belongs to the Hitchin component.

\section{Boundary of the Gothen components}\label{sec:background}

\subsection{Length spectrum of a representation in the Gothen components} 

Given a Higgs bundle parametrized by $(\mu, \nu)$ in a Gothen component $\mathpzc{G}(k)$ with $g-1<k<3g-3$, we consider the Riemannian metric $h$ on $S$ written locally as $e^{\psi_{1}-\psi_{2}}|dz|^{2}$, where the pair $(\psi_{1}, \psi_{2})$ is the solution to the Hitchin equations \eqref{eq:Hitchin}. It is not hard to show (see for instance \cite{TW}) that, up to a multiplicative constant, $h$ is the induced metric on the associated maximal surface in $\h^{2,2}$ found in \cite{BTT}. We denote by $\mathcal{M}(k)$ the space of such metrics up to isotopy. In order for this matching to be well-defined, we need to make sure that $h$ does not depend on the particular pair $(\mu, \nu)$ chosen, as pairs that are in the same $\mathbb{C}^{*}$-orbit give rise to the same representation.

\begin{lemma}\label{lm:well_def} Assume that $(\mu, \nu)$ and $(\widetilde{\mu}, \widetilde{\nu})$ lie in the same $\mathbb{C}^{*}$-orbit. Then $h=\widetilde{h}$.
\end{lemma}
\begin{proof} It is sufficient to show that if $\widetilde{\mu}=\lambda\mu$ and $\widetilde{\nu}=\lambda^{-1}\nu$ for some $\lambda \in \mathbb{C}^{*}$, then $\widetilde{\psi}_{1}=\psi_{1}-\log(|\lambda|)$ and $\widetilde{\psi}_{2}=\psi_{2}-\log(|\lambda|)$ solve the system
\begin{equation*}
\begin{cases}
\Delta \widetilde{\psi}_{1} =e^{\widetilde{\psi}_{1} - \widetilde{\psi}_{2}}-|\widetilde{\nu}|^{2} e^{-2\widetilde{\psi}_{1}}  \\ 
\Delta \widetilde{\psi}_{2} =|\widetilde{\mu}|^{2}e^{2\widetilde{\psi}_{2}}- e^{\widetilde{\psi}_{1} - \widetilde{\psi}_{2}} 
\end{cases} \ .
\end{equation*}
This follows from a simple direct computation.
\end{proof}

The following result shows that $h$ is negatively curved, which is a fundamental step in our construction. The proof is an easy adaptation of the arguments in \cite{OT}.

\begin{prop}\label{prop:neg_curv} Let $\psi_{1}$ and $\psi_{2}$ be the solutions to Equations (\ref{eq:Hitchin}) with given data $(\mu, \nu)$. Then the metric $e^{\psi_{1}-\psi_{2}}|dz|^{2}$ is strictly negatively curved.
\end{prop}
\begin{proof} The curvature of $h$ is
\[
    \kappa(h)=-2\Delta_{h}\log(h)=2\left(\frac{e^{-2\psi_{1}}|\nu|^{2}}{e^{\psi_{1}-\psi_{2}}}+\frac{e^{2\psi_{2}}|\mu|^{2}}{e^{\psi_{1}-\psi_{2}}}-2\right),
\]
hence it is strictly negative if and only if $f_{1}+f_{2}<2$, where 
\[
    f_{1}=\frac{e^{-2\psi_{1}}|\nu|^{2}}{e^{\psi_{1}-\psi_{2}}} \ \ \ \ \text{and} \ \ \ \ f_{2}=\frac{e^{2\psi_{2}}|\mu|^{2}}{e^{\psi_{1}-\psi_{2}}} \ .
\]
Now, outside the zeros of $\mu$ and $\nu$ we have that
\begin{align*}
    \Delta_{h}\log(f_{1})&=-3\Delta_{h}\psi_{1}+\Delta_{h}\psi_{2}=3\frac{e^{-2\psi_{1}}|\nu|^{2}}{e^{\psi_{1}-\psi_{2}}}+\frac{e^{2\psi_{2}}|\mu|^{2}}{e^{\psi_{1}-\psi_{2}}}-4=-4+3f_{1}+f_{2} \\
    \Delta_{h}\log(f_{2})&=3\Delta_{h}\psi_{2}-\Delta_{h}\psi_{1}=3\frac{e^{2\psi_{2}}|\mu|^{2}}{e^{\psi_{1}-\psi_{2}}}+\frac{e^{-2\psi_{1}}|\nu|^{2}}{e^{\psi_{1}-\psi_{2}}}-2=3f_{2}+f_{1}-2 \ .
\end{align*}
Then
\begin{align*}
    \Delta_{h}\log(f_{1}+f_{2})&\geq \frac{f_{1}\Delta_{h}\log(f_{1})+f_{2}\Delta_{h}\log(f_{2})}{f_{1}+f_{2}}\\
    &=\frac{-4f_{1}+3f_{1}^{2}+2f_{1}f_{2}+3f_{2}^{2}+f_{1}f_{2}-2f_{2}}{f_{1}+f_{2}}\\
    &>\frac{-4f_{1}-4f_{2}+3f_{1}^{2}+2f_{1}f_{2}+3f_{2}^{2}+f_{1}f_{2}}{f_{1}+f_{2}}\geq -4+2(f_{1}+f_{2}) \ .
\end{align*}
We deduce that, if $f_{1}+f_{2}$ takes its maximum outside the zeros of $\mu$ and $\nu$, then $f_{1}+f_{2}<2$ and $\kappa(h)<0$. On the other hand, it is clear that $f_{1}+f_{2}$ cannot be maximal at a common zero of $\mu$ and $\nu$ because it is always non-negative. Then, if $f_{1}+f_{2}$ takes its maximum at a zero $p$ of $\nu$ but not of $\mu$, we have $f_{1}(p)=0$ and $p$ is a maximum of $f_{2}$. But then,
\[
    0\geq \Delta_{h}\log(f_{2})(p)=-2+3f_{2}(p)+f_{1}(p)=-2+3f_{2}(p)
\]
which implies that $f_{2}(p)\leq \frac{2}{3}$ and $f_{1}+f_{2}\leq f_{1}(p)+f_{2}(p)\leq \frac{2}{3}$ everywhere on $S$. Hence $\kappa(h)<0$ in this case as well. Similarly, if $f_{1}+f_{2}$ takes its maximum at a zero $p$ of $\mu$ but not of $\nu$, then $f_{2}(p)=0$ and $p$ is a maximum of $f_{1}$. But then,
\[
    0\geq \Delta_{h}\log(f_{1})(p)=-4+3f_{1}(p)+f_{2}(p)=-4+3f_{1}(p)
\]
which implies that $f_{1}(p)\leq \frac{4}{3}$ and $f_{1}+f_{2}\leq f_{1}(p)+f_{2}(p)\leq \frac{4}{3}$ and the curvature is negative everywhere on $S$. 
\end{proof}

Negative curvature guarantees that for every $\gamma \in \pi_{1}(S)$ there is a unique geodesic representative for $h$ in its free homotopy class. We denote by $\ell_{h}(\gamma)$ its length and define the length spectrum of the Gothen representation $\rho \in \mathpzc{G}(k)$ corresponding to the Higgs bundle with parameter $(\mu, \nu)$ as the collection $\{\ell_{h}(\gamma)\}_{\gamma \in \pi_{1}(S)}$. \\
By a result of Otal (\cite{Otal}), the length spectrum of $h$ can be realized by a geodesic current $L_{h}$, thus we can map each Gothen component $\mathpzc{G}(k)$ inside the space of geodesic currents $\mathcal{C}(S)$ (see \cite{bonahonbouts}, \cite{Bonahon_currents}). In Section \ref{subsec:limits} we will prove that this inclusion is proper and describe its projective closure.

\subsection{Comparison results} We collect here some results about how the metric $h=e^{\psi_{1}-\psi_{2}}|dz|^{2}$ compares with the quartic differential metric $|\mu\nu|^{\frac{1}{2}}$ and the hyperbolic metric in the conformal class. Most of the results here already appeared in \cite{QL_cyclic} in some form.

\begin{lemma}\label{lm:lower_bound_flat} Let $h=e^{\psi_{1}-\psi_{2}}|dz|^{2}$ be the metric associated to the pair $(\mu, \nu)$. Let $q=\mu\nu$ be the holomorphic quartic differential in the Hitchin base. Then $|q|^{\frac{1}{2}}\leq h$.
\end{lemma}
\begin{proof}
We consider the function $u=e^{4\psi_{2}-4\psi_{1}}|q|^{2}$, which is well-defined everywhere on $S$. It is clear that $u\geq 0$ and takes a maximum outside the zeros of $q$. Now, outside the zeros of $q$, the function $u$ satisfies
\begin{align*}
    \Delta \log(u) &= 4\Delta \psi_{2}-4\Delta \psi_{1} \\
    &= 4|\mu|^{2}e^{2\psi_{2}}-8e^{\psi_{1}-\psi_{2}}+4|\nu|^{2}e^{-2\psi_{1}} \\
    &\geq 8|q|e^{\psi_{2}-\psi_{1}}-8e^{\psi_{1}-\psi_{2}} \\
    &= 8e^{\psi_{1}-\psi_{2}}(e^{2\psi_{2}-2\psi_{1}}|q|-1)\\
    &= 8e^{\psi_{1}-\psi_{2}}(\sqrt{u}-1) .
\end{align*}
Hence, at a maximum of $u$, we have
\[
    0\geq \Delta_{h}\log(u) \geq 8e^{\psi_{1}-\psi_{2}}(\sqrt{u}-1) \ .
\]
We deduce that $u\leq 1$, which implies that
\[
    |q|^{\frac{1}{2}}\leq e^{\psi_{1}-\psi_{2}} |dz|^{2}=h \ .
\]
\end{proof}

\begin{lemma}\label{lm:compare_Hitchin} Let $h=e^{\psi_{1}-\psi_{2}}|dz|^{2}$ and $\widetilde{h}=e^{\widetilde{\psi}_{1}-\widetilde{\psi}_{2}}|dz|^{2}$ be the metrics associated to the pairs  $(\mu, \nu)$ and $(1,q)$ with $q=\mu\nu$, respectively. Then $h\leq \widetilde{h}$.
\end{lemma}
\begin{proof} The function $u=\psi_{1}-\psi_{2}-\widetilde{\psi}_{1}+\widetilde{\psi}_{2}$ satisfies the equation
\begin{align*}
    \Delta u &= \Delta \psi_{1}-\Delta \psi_{2} -\Delta\widetilde{\psi}_{1}+ \Delta \widetilde{\psi}_{2} \\
    &= 2e^{\psi_{1}-\psi_{2}}-|\nu|^{2}e^{-2\psi_{1}}-|\mu|^{2}e^{2\psi_{2}}-2e^{\widetilde{\psi}_{1}-\widetilde{\psi}_{2}}+|q|^{2}e^{-2\widetilde{\psi}_{1}}+e^{2\widetilde{\psi}_{2}} \\
    &=2e^{\widetilde{\psi}_{1}-\widetilde{\psi}_{2}}(e^{u}-1)-(|\nu|e^{-\psi_{1}}-|\mu|e^{\psi_{2}})^{2}+(|q|e^{-\widetilde{\psi}_{1}}-e^{\widetilde{\psi}_{2}})^{2}+2|q|e^{-\widetilde{\psi}_{1}+\widetilde{\psi}_{2}}(1-e^{-u}) \\
    &\geq 2e^{\widetilde{\psi}_{1}-\widetilde{\psi}_{2}}(e^{u}-1) + 2|q|e^{-\widetilde{\psi}_{1}+\widetilde{\psi}_{2}}(1-e^{-u})
\end{align*}
because by the proof of \cite[Theorem 6.2, page 1257]{QL_cyclic} the term $(|q|e^{-\widetilde{\psi}_{1}}-e^{\widetilde{\psi}_{2}})^{2}-(|\nu|e^{-\psi_{1}}-|\mu|e^{\psi_{2}})^{2}$ is positive. It follows that the constant function $0$ is a supersolution, so $u\leq 0$ and $h=e^{\psi_{1}-\psi_{2}}|dz|^{2}=e^{u}e^{\widetilde{\psi}_{1}-\widetilde{\psi}_{2}}|dz|^{2}\leq \widetilde{h}$.
\end{proof}

\begin{rmk}The metric $\tilde{h}$ is, up to a multiple, the induced metric on the maximal surface in $\h^{2,2}$ that is equivariant under a Hitchin representation. This metric was studied extensively in \cite{OT_Sp4}.
\end{rmk}

\begin{cor}\label{cor:area_bound} The area of $S$ endowed with the metric $h$ satisfies
\[
    \Area(S, |q|^{\frac{1}{2}}) \leq \Area(S,h) \leq \frac{3}{2}\Area(S, |q|^{\frac{1}{2}})+\frac{3\pi}{2}|\chi(S)| \ .
\]
\end{cor}
\begin{proof} The lower bound follows immediately from Lemma \ref{lm:lower_bound_flat}. For the upper bound we use Lemma \ref{lm:compare_Hitchin} and \cite[Proposition 3.9]{OT_Sp4}.
\end{proof}

\begin{lemma}\label{lm:lower_bound_hyp} Let $\sigma=\sigma(z)|dz|^{2}$ be the hyperbolic metric in the conformal class. Then $\frac{1}{4}\sigma \leq h$.
\end{lemma}
\begin{proof}
Consider the function $f=e^{\psi_{1}-\psi_{2}}\sigma^{-1}$, which is well-defined everywhere on $S$. Then,
\begin{align*}
    \Delta \log(f) &=\Delta \psi_{1}-\Delta \psi_{2}-\Delta \log(\sigma)\\
    &= -e^{2\psi_{2}}|\mu|^{2}+2f\sigma-|\nu|^{2}e^{-2\psi_{1}}-\frac{1}{2}\sigma  \\ &\leq \left(2f-\frac{1}{2}\right)\sigma 
\end{align*}
By the maximum principle, $f\geq \frac{1}{4}$ and $h=f\sigma\geq \frac{1}{4} \sigma$
\end{proof}

\subsection{Projective limits in the space of currents}\label{subsec:limits}
By Proposition \ref{prop:neg_curv}, the metric $h$ is negatively curved, thus we can embed $\mathcal{M}(k)$ into the space of geodesic currents for all $g-1<k<3g-3$. We will denote by $L_{h}$ the associated geodesic current. We recall the two main features of this current:
\begin{enumerate}[i)]
    \item for every curve $\gamma \in \pi_{1}(S)$, we have $\ell_{h}(\gamma)=i(L_{h},\delta_{\gamma})$, where $\ell_{h}(\gamma)$ denotes the length of the unique geodesic representative of $\gamma$ for $h$;
    \item $i(L_{h},L_{h})=\frac{\pi}{2}\Area(S,h)$ \ .
\end{enumerate}

\begin{thm} \label{thm:limit_induced} Let $\rho_{n}\in \mathpzc{G}(k)$ be a sequence of Gothen representations leaving every compact set. Let $(\mu_{n}, \nu_{n})$ be the sequence parameterizing the corresponding Higgs bundles and let $h_{n}$ be the corresponding sequence of Riemannian metrics on $S$. Then there is a sequence of positive real numbers $t_{n}$ and a mixed structure $\eta$ so that $t_{n}L_{h_{n}}\to \eta$. 
\end{thm}
\begin{proof} Set $q_{n}=\mu_{n}\nu_{n}$. We distinguish two cases, depending on whether the area $\|q_{n}\|$ is uniformly bounded or not. \\
\underline{\textit{First case:} $\sup\|q_{n}\|<\infty$.} By Corollary \ref{cor:area_bound}, the self-intersection $i(L_{h_{n}},L_{h_{n}})$ is uniformly bounded. We notice that, because $\rho_{n}$ leaves every compact set, the sequence of hyperbolic metrics $\sigma_{n}$ in the conformal class of $h_{n}$ must necessarily diverge. Otherwise, up to subsequences, we could assume that $\sigma_{n}\to \sigma \in \mathcal{T}(S)$ and we could write $q_{n}=\tau_{n}\tilde{q}_{n}$, where 
\[
    \tau_{n}=\|q_{n}\|_{\infty}:=\max_{S}\frac{|q_{n}|^{2}}{\sigma_{n}^{4}}
\]
and $\tilde{q}_{n}$ converges uniformly to a non-vanishing quartic differential $\tilde{q}_{\infty} \in \mathcal{Q}(S,\sigma_{\infty})$, as unit spheres of holomorphic differentials over the thick part of Teichm\"uller space are compact. Since
\[
    \|q_{n}\|=\int_{S}|\tau_{n}|^{\frac{1}{2}}|\tilde{q}_{n}|^{\frac{1}{2}}=|\tau_{n}|^{\frac{1}{2}}\int_{S} \frac{|\tilde{q}_{n}|^{\frac{1}{2}}}{\sigma_{n}} dA_{\sigma_{n}} 
\]
and
\[
    \frac{|\tilde{q}_{n}|^{\frac{1}{2}}}{\sigma_{n}} dA_{\sigma_{n}} \to \frac{|\tilde{q}_{\infty}|^{\frac{1}{2}}}{\sigma_{\infty}} dA_{\sigma_{\infty}} \neq 0 \ ,
\]
the bound on $\|q_{n}\|$ implies that $\tau_{n}$ is uniformly bounded and $q_{n}$ converges up to subsequences. Because $q_{n}=\mu_{n}\nu_{n}$ and, up to the $\mathbb{C}^{*}$-action we can assume that $|\mu_{n}|$ and $|\nu_{n}|$ are uniformly bounded, the pair $(\mu_{n}, \nu_{n})$ would also converge up to subsequence, thus contradicting the fact that the representations $\rho_{n}$ leave every compact set. Therefore, the sequence $\sigma_{n}$ of hyperbolic metrics in the conformal class of $h_{n}$ diverges, as claimed. Then, by Lemma \ref{lm:lower_bound_hyp}, the length spectrum of $L_{h_{n}}$ is unbounded. Since $\Pp\mathcal{C}(S)$ is compact, there exists a sequence $t_{n}\to 0$ such that $t_{n}L_{h_{n}}\to L_{\infty}$. We easily deduce that $i(L_{\infty},L_{\infty})=0$, hence $L_{\infty}$ is a measured lamination, which we can interpret as a mixed structure with no flat parts. \\
\underline{\textit{Second case:} $\sup\|q_{n}\|=\infty$.} By Corollary \ref{cor:area_bound}, the self-intersection of $L_{h_{n}}$ diverges as $\|q_{n}\|$. Since every geodesic current has finite self-intersection, we need to rescale $L_{h_{n}}$ at least by $\frac{1}{\sqrt{\|q_{n}\|}}$. Let us denote
\[
    \hat{L}_{h_{n}}=\frac{1}{\sqrt{\|q_{n}\|}}L_{h_{n}} \ .
\]
If the length spectrum of $\hat{L}_{h_{n}}$ is still unbounded, then there is a sequence $t_{n}\to 0$ such that $t_{n}\hat{L}_{h_{n}} \to \hat{L}_{\infty}$, which has vanishing self-intersection, thus $L_{h_{n}}$ converges projectively to a measured lamination.

If the length spectrum of $\hat{L}_{h_{n}}$ is uniformly bounded, then by Lemma \ref{lm:lower_bound_flat}, the length spectrum of the unit area flat metrics $|q_{n}|^{\frac{1}{2}}/\|q_{n}\|$ is uniformly bounded as well. Thus, by \cite[Theorem 2.5]{OT_Sp4}, the geodesic currents $L_{q_{n}}$ converges in $\Pp\mathcal{C}(S)$ to a mixed structure $m$ that is not purely laminar. This furnishes an orthogonal decomposition (for the intersection form $i$)  of the surface $S$ into a collection of incompressible subsurfaces $\{S_{j}'\}_{j=1}^{l}$, obtained by cutting $S$ along the simple closed curves $\gamma_{i}$, for which $m$ is induced by a flat metric on each $S_{j}'$ and is a measure lamination on the complement. Moreover, we can assume that each simple closed curve $\gamma_{i}$ bounds at least one flat part, induced by a meromorphic quartic differential $\tilde{q}_{j}$. On each $S_{j}'$, the inequalities $|q_{n}|^{\frac{1}{2}}\leq h_{n} \leq \widetilde{h}_{n}$ given by Lemma \ref{lm:compare_Hitchin} and Lemma \ref{lm:lower_bound_flat}, together with \cite[Corollary 3.14]{OT_Sp4}, imply the sequence $h_{n}$, renormalized by $\|q_{n}\|$, converges to $|\tilde{q}_{j}|^{\frac{1}{2}}$ uniformly on compact sets outside the zeros and poles of $\tilde{q}_{j}$. We deduce that on each $S_{j}'$, we have
\[
    \hat{L}_{\infty}=\lim_{n \to +\infty}\hat{L}_{h_{n}}=\lim_{n\to +\infty}\frac{1}{\sqrt{\|q_{n}\|}}L_{q_{n}}=L_{\tilde{q}_{j}} \ ,
\]
because uniform convergence of metrics implies convergence in the length spectrum (\cite[Proposition 5.3]{Charles_dPSL}). In particular, we notice that
\[
    \lim_{n \to +\infty}i(\hat{L}_{h_{n}}, \delta_{\gamma_{i}})=0
\]
so that the same collection of curves $\gamma_{i}$ can be used for the orthogonal decomposition of $\hat{L}_{\infty}$. Therefore, we can write
\[
    \hat{L}_{\infty}=\sum_{j=1}^{l}L_{\tilde{q}_{j}}+\lambda \ ,
\]
where $\lambda$ is a geodesic current supported in the complement of $\bigcup_{j}S_{j}'$. We claim that $\lambda$ is a measure lamination: in fact using again the inequalities $|q_{n}|^{\frac{1}{2}} \leq h_{n}\leq \tilde{h}_{n}$ and \cite[Lemma 3.15]{OT_Sp4}, we have
\begin{align*}
    \frac{\pi}{2}&=\frac{\pi}{2}\lim_{n\to +\infty}\frac{\Area(S,h_{n})}{\|q_{n}\|}
    =\lim_{n \to +\infty}i(\hat{L}_{h_{n}},\hat{L}_{h_{n}})\\
    &=i(\hat{L}_{\infty}, \hat{L}_{\infty})
    =\sum_{j=1}^{l} i(L_{\tilde{q}_{j}}, L_{\tilde{q}_{j}})+i(\lambda, \lambda)\\
    &=\lim_{n\to +\infty} \frac{1}{\|q_{n}\|}i(L_{q_{n}},L_{q_{n}})+i(\lambda, \lambda)=\frac{\pi}{2}+i(\lambda, \lambda) \ ,
\end{align*}
so $\lambda$ has vanishing self-intersection. 
\end{proof}

\begin{proof}[Proof of Theorem A] By Theorem \ref{thm:limit_induced}, we know that $\partial \overline{\mathcal{M}(k)}\subseteq \Pp\Mix_{4}(S)$. Consider now the family of metrics $h_{t}$ associated to a ray $tq$ in the Hitchin base for a fixed unit area quartic differential. By the proof of Theorem \ref{thm:limit_induced}, we know that $L_{h_{t}}$ converges projectively to $L_{q}$ as $t\to +\infty$. Therefore, $\partial \overline{\mathcal{M}(k)}\supseteq \overline{\Flat_{4}(S)}=\Pp\Mix_{4}(S)$ for all $g-1<k<3g-3$. Since the map $\mathpzc{G}(k) \rightarrow \mathpzc{M}(k)$ is proper by Theorem \ref{thm:limit_induced}, we can identify the boundary of the Gothen components $\mathpzc{G}(k)$ with $\mathbb{P}\mathrm{Mix}_{4}(S)$.
\end{proof}

\section{Induced metrics on the minimal surfaces}
Given a maximal representation $\rho:\pi_{1}(S) \rightarrow \Sp(4,\R)$, there is a unique $\rho$-equivariant conformal harmonic map $f_{\rho}:\tilde{S} \rightarrow \Sp(4,\R)/\mathrm{U}(2)$ (\cite{Corlette}, \cite{collier2016maximal}). The image of $f_{\rho}$ is thus a minimal surface in the symmetric space and one could have defined the length spectrum of the maximal representation by considering the induced metric on this surface. In this section we show that how this would have still led to identifying the boundary of the Gothen and Hitchin components with the space of projectivized mixed structures. \\

Dai and Li (\cite{QL_cyclic}) write an explicit expression for the induced metric $g$ on the minimal surface in terms of the Higgs bundle data: in local coordinates
\[
    g=16(|\nu|^{2}e^{-2\psi_{1}}+2e^{\psi_{1}-\psi_{2}}+|\mu|^{2}e^{2\psi_{2}})|dz|^{2}
\]
where $\psi_{1}$ and $\psi_{2}$ are the solutions to Hitchin equations (\ref{eq:Hitchin}). We denote by $\tilde{g}$ the induced metric on the minimal surface associated to the Higgs bundle in the Hitchin component in the same fiber, i.e.
\[  
    \tilde{g}=16(|q|^{2}e^{-2\tilde{\psi}_{1}}+2e^{\tilde{\psi}_{1}-\tilde{\psi}_{2}}+e^{2\tilde{\psi}_{2}})|dz|^{2}
\]
where $q=\mu\nu$.

Let us start by analyzing the behavior of the metric $\tilde{g}$. By \cite[Theorem 5.6]{QL_cyclic}, the metric $\tilde{g}$ is strictly negatively curved and thus there is a unique geodesic current $L_{\tilde{g}}$ that records its marked length spectrum. 

\begin{thm}\label{thm:B_Hitchin} Let $\rho_{n}$ be a sequence of representations leaving every compact set in the Hitchin component and let $L_{\tilde{g}_{n}}$ be the corresponding sequence of geodesic currents. Then there is a sequence $t_{n} \to 0$ and $\eta \in \mathbb{P}\Mix(S)$ such that $t_{n}L_{\tilde{g}_{n}}$ converges to $\eta$.
\end{thm}
\begin{proof} We first observe that, by the proof of Proposition \ref{prop:neg_curv}, we know that 
\[
    |q|^{2}e^{-2\tilde{\psi}_{1}}+e^{2\tilde{\psi}_{2}} \leq 2e^{\tilde{\psi}_{1}-\tilde{\psi}_{2}} \ ,
\]
hence $\tilde{g}$ is uniformly bi-Lipschitz to $\tilde{h}$, as \begin{equation}\label{eq:bi_Lip}
32\tilde{h} \leq \tilde{g} \leq 64\tilde{h}.
\end{equation}
In particular, the length spectrum of $\tilde{g}_{n}$ diverges at the same rate as the length spectrum of $\tilde{h}_{n}$. Let $t_{n}$ be the scaling factors such that the geodesic currents $L_{\tilde{h}_{n}}$ converge to a projectivized mixed structure $\eta'$. Then, because the length spectrum of $t_{n} L_{\tilde{g}_{n}}$ is uniformly bounded, up to subsequences, $t_{n}L_{\tilde{g}_{n}}$ converges to a geodesic current $L_{\infty}$. Now, if $\eta'$ is purely laminar, then 
\[
    0=\iota(\eta', \eta')=\lim_{n \to +\infty}t_{n}^{2}\iota(L_{\tilde{h}_{n}}, L_{\tilde{h}_{n}})=\lim_{n \to +\infty} \frac{t_{n}^{2}\pi}{2}\Area(S, \tilde{h}_{n})
\]
and by Equation (\ref{eq:bi_Lip}) we deduce that
\[
    0=\lim_{n \to +\infty}\frac{t_{n}^{2}\pi}{2}\Area(S, \tilde{g}_{n})=\iota(L_{\infty}, L_{\infty})
\]
so $L_{\infty}$ is a measured lamination as well. 
Assume now that $\eta'$ is not purely laminar and let $S_{j}$ be a subsurface on which $\eta'$ is given by a meromorphic quartic differential of finite area $q_{j}$. Note that this case can happen (see \cite[Theorem 3.16]{OT_Sp4}) only when the area of the quartic differentials metrics $|q_{n}|^{\frac{1}{4}}$ is unbounded, and such $q_{j}$ arises then as limit of the rescaled flat metrics $t_{n}^{2}|q_{n}|^{\frac{1}{2}}$. Set $\tilde{u}_{1,n}=\tilde{\psi}_{1,n}-\log(|q_n|^{\frac{3}{4}})$ and $\tilde{u}_{2,n}=\tilde{\psi}_{2,n}-\log(|q_n|^{\frac{1}{4}})$. We can then re-write $\tilde{g}_n$ as
\[
    \tilde{g}_{n}=16|q_{n}|^{\frac{1}{2}}(e^{-2\tilde{u}_{1,n}}+2e^{\tilde{u}_{1,n}-\tilde{u}_{2,n}}+e^{2\tilde{u}_{2,n}})
\]
By \cite[Corollary 3.13]{OT_Sp4}, the sequences $u_{1,n}$ and $u_{2,n}$ converge to $0$ uniformly on compact sets outside the zeros and poles of $q_{j}$, whereas the renormalized flat metrics $t_{n}^{2}|q_{n}|^{\frac{1}{2}}$ converge to $|q_{j}|^{\frac{1}{2}}$ by assumption. Hence, we can conclude that on $S_{j}$ we have
\[
    L_{{\infty}_{|_{S_{j}}}}=8L_{q_{j}}=8 {\eta'}_{|_{S_{j}}} \ . 
\]
Therefore, we can write $L_{\infty}=\lambda+\sum_{j}8L_{q_{j}}$, where $\lambda$ is a geodesic current supported on the complement of the flat subsurfaces $\bigcup S_{j}$ in the decomposition of $\eta'$. We claim that $\lambda$ is a measured lamination. Indeed, recalling that $\iota(\eta', \eta')=\frac{\pi}{2}$ (see \cite[Thoerem 3.16]{OT_Sp4}), we have
\begin{align*}
    \frac{64\pi}{2} &=\iota(8\eta', 8\eta')=64\lim_{n \to +\infty}t_{n}^{2}\iota(L_{\tilde{h}_{n}}, L_{\tilde{h}_{n}})=64\lim_{n\to +\infty}\frac{t_{n}^{2}\pi}{2}\Area(S, \tilde{h}_{n}) \\
    & \geq \lim_{n \to +\infty} \frac{t_{n}^{2}\pi}{2}\Area(S, \tilde{g}_{n})=\lim_{n \to +\infty}t_{n}^{2}\iota(L_{\tilde{g}_{n}}, L_{\tilde{g}_{n}})= \iota(L_{\infty}, L_{\infty})\\
   & =\frac{64\pi}{2}+\iota(\lambda, \lambda).
\end{align*}
We conclude that $\lambda$ has vanishing self-intersection and thus is a measured lamination. If we further renormalize $L_{\infty}$ by dividing by $8$, we obtain a mixed structure $\eta$ with self-intersection $\frac{\pi}{2}$, in other words, the areas of the flat subsurfaces sum up to $1$.
\end{proof}

We now deduce Theorem \ref{thmB} for induced metrics arising from representations in the Gothen components.

\begin{proof}[Proof of Theorem \ref{thmB}]
Since Equation (\ref{eq:bi_Lip}) holds for the Gothen components as well, as it is based on Proposition \ref{prop:neg_curv}, we know that the length spectrum of $g_{n}$ grows at the same rate as the length spectrum of the metrics $h_{n}$, thus if we choose scaling factors such that $t_{n}L_{h_{n}}$ converges to a mixed structure $\eta$, then $t_{n}L_{g_{n}}$  must converge as well. Let $L_{\infty}$ be such limit. The same argument as in the proof of Theorem \ref{thm:B_Hitchin} shows that if $\eta$ is purely laminar, then $L_{\infty}$ is a measured lamination. 

Let now $S_{j}$ be a subsurface in which $\eta$ is the Liouville current of a flat metric induced by a meromorphic quartic differential $q_{j}$. By the proof of \cite[Theorem 6.2, page 1257]{QL_cyclic}, the following inequalities hold
\[
    16(|q_{n}|^{2}e^{-2\tilde{\psi}_{1}}+ 2|q_{n}|^{\frac{1}{2}}+|q_{n}|^{2}e^{-2\tilde{\psi}_{1}}) \leq g_{n} \leq \tilde{g}_{n} \ .
\]
Introducing again the function $\tilde{u}_{1,n}=\tilde{\psi}_{1,n}-\log(|q_n|^{\frac{3}{4}})$, we can rewrite the left-hand side as
\[
    32|q_{n}|^{\frac{1}{2}}(e^{-2\tilde{u}_{1,n}}+1) \ .
\]
Because $\tilde{u}_{1,n}$ converges uniformly to $0$ outside the zeros and poles of $q_{j}$ and the metrics $\tilde{g_{n}}$ and $|q_{n}|^{\frac{1}{2}}$ rescaled by $t_{n}^{2}$ both converge to $64|q_{j}|^{\frac{1}{2}}$ and $|q_{j}|^{\frac{1}{2}}$ respectively, we can conclude that $g_{n}$ converges to $64|q_{j}|^{\frac{1}{2}}$ and thus $L_{\infty}$ restricted to the subsurface $S_{j}$ coincides with $64L_{q_{j}}$. If we then write $L_{\infty}=\sum L_{q_{j}} + \lambda$, where $\lambda$ is a measured lamination supported on on the complement of the subsurfaces $S_{j}$, then the same area argument as in the proof of Theorem \ref{thm:B_Hitchin} shows that $\lambda$ has self-intersection $0$ and thus is a measured lamination. Hence, $\eta=\frac{1}{64}L_{\infty}$ is a mixed structure where the flat pieces have unit area.
\end{proof}

\bibliographystyle{alpha}
\bibliography{bs-bibliography}

\bigskip
\noindent \footnotesize \textsc{Department of Mathematics and Statistics, University of Massachusetts, Amherst}\\
\emph{E-mail address:}  \verb|ouyang@math.umass.edu|

\bigskip
\noindent \footnotesize \textsc{Department of Mathematics, Rice University}\\
\emph{E-mail address:} \verb|andrea_tamburelli@libero.it|

\end{document}